\documentclass [a4paper,10pt]{article}
\usepackage [utf8] {inputenc}
\usepackage {geometry}
\usepackage {amsmath}
\usepackage {amsthm}
\usepackage {amsfonts}
\title {A variance for $k$-free numbers in arithmetic progressions of given modulus}

\date {Tomos Parry}
\begin {document}
\maketitle
\begin {center}
$\hspace {1mm}$
\\
\section {- \hspace {2mm}Introduction}
\end {center}
Let
\[ \mathcal S=\{ n\in \mathbb N|\text { there is no prime $p$ with }p^k|n\} ,\]
the set of $k$-free numbers.  For 
some suitable main term $\eta (q,a)$ to be defined soon enough we will study in this paper the object
\[ \sum _{a=1}^q\left (\sum _{n\leq x\atop {n\in \mathcal S\atop {n\equiv a\text { mod }(q)}}}1-x\eta (q,a)\right )^2,\]
a variance for $k$-free numbers in arithmetic progressions when averaging over a (complete) residue system.
One would like to establish for some $q$ that this is
\[ \approx q\left (\frac {x}{q}\right )^{1/k}\]
since this would mean that on average
\begin {equation}\label {conj}
\sum _{n\leq x\atop {n\in \mathcal S\atop {n\equiv a\text { mod }(q)}}}1-x\eta (q,a)\approx \left (\frac {x}{q}\right )^{1/2k}.
\end {equation}
Since an improvement in the error term in the classical statement
\[ \sum _{n\leq x\atop {n\in \mathcal S}}1=\frac {x}{\zeta (k)}+\mathcal O\left (x^{1/k}\right )\]
is tantamount to a better zero-free region for the zeta function, information as to the true 
size of the LHS of \eqref {conj} is relevant.
\\
\\ Averaging just over the reduced classes an asymptotic formula for the variance, in the squarefree case, is 
already established in \cite {ramon} with error essentially
\begin {equation}\label {banana}
\ll q\left (\frac {x}{q}\right )^{1/3}+\bigg (\frac {x}{q}\bigg )^{23/15}.
\end {equation}
Before this only upper bound results are recorded (see \cite {hecke} and the references therein), although these are stronger 
in the range where the above asymptotic formulas don't hold and are concerned with more general sequences than 
the squarefrees.  In this paper, we improve the first error term in \eqref {banana}.
\\
\newtheorem* {theo}{Theorem}
\begin {theo}\label {theo}
Let $k\geq 2$ and denote by $\mathcal S$ the set of $k$-free numbers.  For $q,a\in \mathbb N$ and $x>0$ define
\begin {eqnarray}\label {ramonramon}
\eta (q,a)=\sum _{d=1\atop {(q,d^k)|a}}^\infty \frac {\mu (d)}{[q,d^k]},\hspace {10mm}
E_x(q,a)=\sum _{n\leq x\atop {n\in \mathcal S\atop {n\equiv a(q)}}}1-x\eta (q,a)\hspace {15mm}
\end {eqnarray}
and 
\begin {equation}\label {ramonramon2}
V_x(q)=\sum _{a=1}^q|E_x(q,a)|^2.
\end {equation}
Define
\[ C_k=\frac {2k}{(1/k-1)\zeta (2)}\prod _p\frac {1-2/(p^k+p^{k-1})}{1-p^{1-1/k}}\]
%\hspace {5mm}\text { and }\hspace {5mm}
and
\[ f_k(q)=C_k\prod _{p|q}\frac {1-2/p^k+(q,p^k)^{1/k-1}/p}{1-2/p^k+1/p}.\]
For $1\leq q\leq x$ we have for every $\epsilon >0$
\[ V_x(q)=q\left (\frac {x}{q}\right )^{1/k}f_k(q)+
\mathcal O_{k,\epsilon }\left (x^\epsilon \left (q\left (\frac {x}{q}\right )^{2/(9-2/k)}
+\frac {x^{1+2/(k+1)}}{q}\right )\right ).\]
\end {theo}
This is an asymptotic formula for $k=2,3,4$.  The relevance of our result is the improvement in 
the first error term, which for $k=2$ seems decently small.  This is obtained by a careful analysis of the 
integrals arising from an application of Perron's formula.  (Our second error term is weaker than in \eqref {banana} but most 
likely can be made to be just as small for the squarefress by arguing, as in that paper, with the square sieve.)
\\
\\ We consider $k\geq 2$ and $q\leq x$ as fixed throughout.  % and assume $x/q\rightarrow \infty $ since the result is clear otherwise.  E
Each time $\epsilon $ appears it is to be understood that it may be taken arbitrarily small at each occurence.  Fix 
some $0<\delta <1/2k$.  All $\ll ,\mathcal O$ constants depend on $\epsilon ,k$ and $\delta $.
\begin {center}
\section {- \hspace {2mm}Lemmas}
\end {center}
For $\mathfrak R\mathfrak e(s)>1$ define
\[ \mathcal F(s)=\sum _{d,d'=1}^\infty \frac {\mu (d)\mu (d')}{[d^k,d'^k][q,(d^k,d'^k)]^{s}}\]
and for $\mathfrak R\mathfrak e(s)\geq -1+\delta $ define
\[ \mathcal F^*(s)=\prod _{p|q}\frac {1+(q,p^k)^s/p^{k(1+s)}}{1+1/p^{k(1+s)}}
\prod _p\left (1-\frac {2}{p^k\left (1+(q,p^k)^s/p^{k(1+s)}\right )}\right ).\]
The first series converges since the summands are bounded by 
\[ \frac {1}{[d^k,d'^k](d^k,d'^k)};\]
for $\mathfrak R\mathfrak e(s)\geq -1+\delta $
\begin {eqnarray}\label {hwndio}
\left |(q,p^k)^s/p^{k(1+s)}\right |&\leq &\left \{ \begin {array}{ll}1/p^k&\text { for }\mathfrak R\mathfrak e(s)\geq 0
\\ 1/p^{k\delta }&\text { for }\mathfrak R\mathfrak e(s)<0\end {array}\right .\notag 
\\ &\ll &1
\end {eqnarray}
and therefore 
\[ 1+(q,p^k)^s/p^{k(1+s)}\geq 1-1/2^{k\delta }\gg 1\]
so that each Euler factor of the infinite product in $\mathcal F^*(s)$ is of the form 
\[ 1+\mathcal O\left (1/p^k\right )\]
and therefore this product converges and is uniformly bounded 
for $\mathfrak R\mathfrak e(s)\geq -1+\delta $; for $\mathfrak R\mathfrak e(s)\geq -1+\delta $ we have 
%\begin {equation}\label {hwndio2}
\[ \frac {1}{p^{k(1+s)}}\geq \frac {1}{p^{k\delta }}\]
and therefore
\[ 1+1/p^{k(1+s)}\geq 1-1/2^{k\delta }\gg 1\]
%\end {equation}
so each factor in 
the finite product in $\mathcal F^*(s)$ is from \eqref {hwndio} uniformly bounded for $\mathfrak R\mathfrak e(s)\geq -1+\delta $, and 
since we have just said the same is true for the infinite product, we conclude 
that $\mathcal F^*(s)\ll q^\epsilon $ for $\mathfrak R\mathfrak e(s)\geq -1+\delta $.
%\[ \mathcal F^*(s)=\prod _p\left (1-\frac {2}{p^k(1+1/p^{k(1+s)})}\right ).\]
%Denote by $\mathcal E$ the set of $s\in \mathbb C$ for which 
%\[ 1-\frac {2}{p^k}+\frac {1}{p^{k(1+s)}}=0\]
%for some prime $p|q$, and define for $s\in \mathbb C\char92 \mathcal E$
%\[ \Delta (s)=\prod _{p|q}\left (\frac {1-2/p^k+(q,p^k)^{s}/p^{k(1+s)}}{1-2/p^k+1/p^{k(1+s)}}\right );\]
%clearly $\Delta (s)$ is holomorphic on $\mathfrak R\mathfrak e(s)>0$.
\\
\newtheorem {ds}{Lemma}[section]
\begin {ds}\label {ds}
If $\mathfrak R\mathfrak e(s)>1$ then
\[ \mathcal F(s)=\frac {\zeta \left (k(s+1)\right )\mathcal F^*(s)}{q^s\zeta \left (2k(s+1)\right )}.\] 
%If $q$ is odd then $\Delta (s)$ is holomorphic 
%Then $\mathcal F^*(s)$ is holomorphic and uniformly bounded for $\mathfrak R\mathfrak e(s)\geq -1+\delta $.
\end {ds}
\begin {proof}
We have
\begin {eqnarray}\label {yreidal}
\sum _{d,d'}\frac {\mu (d)\mu (d')(d^k,d'^k)^{1-s}(q,d^k,d'^k)^{s}}{d^kd'^k}&=&\sum _{N=1}^\infty \frac {1}{N^k}
\sum _{dd'=N}\mu (d)\mu (d')(d^k,d'^k)^{1-s}(q,d^k,d'^k)^{s}\notag 
\\ &=:&\sum _{N=1}^\infty \frac {a_q(N)}{N^k}.
\end {eqnarray}
Clearly $a_q(N)$ is multiplcative and simple calculations show
\begin {eqnarray*}
a_q(p)=-2,
\end {eqnarray*}
\begin {eqnarray*}
a_q(p^2)=p^{k(1-s)}(q,p^k)^{s}
\end {eqnarray*}
and $a_q(p^t)=0$ for $t\geq 3$.  Consequently
\begin {eqnarray}\label {columbs}
\sum _{N=1}^\infty \frac {a_q(N)}{N^k}&=&
\prod _p\left (1-\frac {2}{p^k}+\frac {(q,p^k)^{s}}{p^{k(1+s)}}\right )\notag 
\\ &=&\prod _p\left (1+\frac {(q,p^k)^{s}}{p^{k(1+s)}}\right )
\prod _p\left (1-\frac {2}{p^k\left (1+(q,p^k)^s/p^{k(1+s)}\right )}\right )\notag 
\\ &=&\prod _p\left (1+\frac {1}{p^{k(1+s)}}\right )\prod _{p|q}\frac {1+(q,p^k)^s/p^{k(1+s)}}{1+1/p^{k(1+s)}}
\prod _p\left (1-\frac {2}{p^k\left (1+(q,p^k)^s/p^{k(1+s)}\right )}\right )\notag 
\\ &=&\frac {\zeta \left (k(1+s)\right )\mathcal F^*(s)}{\zeta \left (2k(1+s)\right )}
%\\ &=&\prod _{p|q}\left (\frac {1-2/p^k+(q,p^k)^{s}/p^{k(1+s)}}{1-2/p^k+1/p^{k(1+s)}}\right )
%\prod _{p}\left (1-\frac {2}{p^k}+\frac {1}{p^{k(1+s)}}\right )\notag 
%\\ &=&\Delta (s)\prod _p\left (1+\frac {1}{p^{k(1+s)}}\right )\prod _p\left (1-\frac {2}{p^k(1+1/p^{k(1+s)})}\right )\notag 
%\\ &=&\frac {\Delta (s)\zeta \left (k(1+s)\right )\mathcal F^*(s)}{\zeta \left (2k(1+s)\right )}
\end {eqnarray}
so that \eqref {yreidal} becomes
\[ \sum _{d,d'}\frac {\mu (d)\mu (d')(d^k,d'^k)^{1-s}(q,d^k,d'^k)^{s}}{d^kd'^k}
=\frac {\zeta \left (k(1+s)\right )\mathcal F^*(s)}{\zeta \left (2k(1+s)\right )}\]
and the claim follows.%  It was already mentioned before the lemma that each Euler product in the infinite product 
%of $\mathcal F^*(s)$ is of the form
%\[ 1+\mathcal O\left (1/p^k\right )\]
%therefore this is obviously holomorphic uniformly bounded.  Since
%\[ \left |\frac {1}{p^{k(1+s)}}\right |\leq \frac {1}{p^{\delta ^k}}\]
%the same is true of the finite product and thereofore of $\mathcal F^*(s)$.
%\\ 
%\\ If $p\geq 3$ and $\mathfrak R\mathfrak e(s)\geq -1+1/2k$ then
%\[ |(1-2/p^k)p^{k(1+s)}|\geq p^{1/2}-2p^{1/2-k}\geq \sqrt 3-\frac {2}{3^{3/2}}>\sqrt 2-2/5>1.01\]
%so that
%\[ 1-\frac {2}{p^k}+\frac {1}{p^{k(1+s)}}\]
%is bounded away from zero.  Consequently $\Delta (s)$ is holomorphic if $q$ only has prime factors $\geq 3$.
\end {proof}
\newtheorem {ds2}[ds]{Lemma}
\begin {ds2}\label {ds2}
%Take $\mathcal F^*(s)$ as in Lemma \ref {ds}.
%(A) For each $n\in \mathbb N$ there are $\lambda _n,W_n\in \mathbb R$ such that
%\[ \mathcal F^*(s)=\sum _{n=1}^\infty \lambda _nW_n^{1+s}\]
%for $\mathfrak R\mathfrak e(s)\geq -1+1/2k$.  Moreover for these $s\leq 0$
%\[ \sum _{n=1}^\infty \left |\lambda _nW_n^{1+s}\right |\ll 1.\]
%(B) Suppose $q$ is odd and has $\omega $ distinct prime 
%factors $p_1,...,p_\omega $.  For integers $l_1,..,l_\omega ,l_1',...,l_\omega '\geq 0$ there 
%are $C_{\mathbf l,\mathbf l'},Z_{\mathbf l,\mathbf l'}\in \mathbb R$ such that
%\[ \Delta (s)=\sum _{l_1,...,l_\omega \geq 0\atop {l_1',...,l_\omega '\geq 0}}
%C_{\mathbf l,\mathbf l'}Z_{\mathbf l,\mathbf l'}^s\]
%for $\mathfrak R\mathfrak e(s)\geq -1+1/2k$.  Moreover for these $s$
%\[ \sum _{l_1,...,l_\omega \geq 0\atop {l_1',...,l_\omega '\geq 0}}\left |
%C_{\mathbf l,\mathbf l'}Z_{\mathbf l,\mathbf l'}^s\right |\ll 1.\]
Suppose $q$ has $\omega $ distinct prime factors $p_1,...,p_\omega $ and let $\mathcal F^*(s)$ be as given at the 
start of this section.  Then for 
each $n\in \mathbb N$ and each $l_1,..,l_\omega ,l_1',...,l_\omega '\geq 0$ there 
are $\lambda _n,W_n,C_{\mathbf l,\mathbf l'},Z_{\mathbf l,\mathbf l'}\in \mathbb R$ with $W_n,Z_{\mathbf l,\mathbf l'}>0$ such that
\[ \mathcal F^*(s)=\sum _{l_1,...,l_\omega \geq 0\atop {l_1',...,l_\omega '\geq 0}}
\sum _{n=1}^\infty C_{\mathbf l,\mathbf l'}Z_{\mathbf l,\mathbf l'}^s\lambda _nW_n^{1+s}\]
for $\mathfrak R\mathfrak e(s)\geq -1+\delta $.  Moreover for $-1+\delta \leq \mathfrak R\mathfrak e(s)\leq 0$ 
\[  \sum _{l_1,...,l_\omega \geq 0\atop {l_1',...,l_\omega '\geq 0}}\sum _{n=1}^\infty 
\left |C_{\mathbf l,\mathbf l'}Z_{\mathbf l,\mathbf l'}^s\lambda _nW_n^{1+s}\right |\ll \log (q+1).\]
\end {ds2}
\begin {proof}
From \eqref {hwndio} we have $|(q,p^k)^s/p^{k(1+s)}|<1$ and therefore
\begin {eqnarray}\label {now}
\prod _p\left (1-\frac {2}{p^k\left (1+(q,p^k)^s/p^{k(1+s)}\right )}\right )
&=&\prod _p\left (1-\frac {2}{p^{k}}\sum _{t\geq 1}\left (\frac {-(q,p^k)^s}{p^{k(1+s)}}\right )^{t-1}\right )\notag 
\\ &=&\sum _{n=1}^\infty f_s^*(n)
\end {eqnarray}
where $f^*_s(n)$ is the multiplicative function given on prime powers by
\[ f_s^*(p^t)=-\frac {2}{p^k}\left (\frac {-(q,p^k)^s}{p^{k(1+s)}}\right )^{t-1}.\]
For any $n\in \mathbb N$ and prime $p|n$ define $t=t(p)$ through $p^t||n$.  Then
\begin {eqnarray}\label {vodan}
f^*_s(n)&=&\prod _{p|n}\left (-\frac {2}{p^k}\right )\left (\frac {-(q,p^k)^s}{p^{k(1+s)}}\right )^{t-1}\notag 
\\ &=&\left (\prod _{p|n}(-1)^{t-1}\right )\left (\prod _{p|n}\frac {-2}{p^{k}}\right )
\left (\prod _{p|n}(q,p^k)^{-(t-1)}\right )
\left (\prod _{p|n}\frac {(q,p^k)^{(t-1)(1+s)}}{p^{(t-1)k(1+s)}}\right ).\hspace {15mm}
\end {eqnarray}
If we now define
\[ \lambda _n=\left (\prod _{p|n}(-1)^{t-1}\right )\left (\prod _{p|n}\frac {-2}{p^{k}}\right )
\left (\prod _{p|n}(q,p^k)^{1-t}\right )\]
and
\[ W_n=\prod _{p|n}\frac {(q,p^k)^{t-1}}{p^{(t-1)k}}\]
then \eqref {vodan} becomes 
\[ f^*(n)=\lambda _nW_n^{1+s}\]
so \eqref {now} becomes
\begin {equation}\label {bisged}
\prod _p\left (1-\frac {2}{p^k\left (1+(q,p^k)^s/p^{k(1+s)}\right )}\right )=\sum _{n=1}^\infty \lambda _nW_n^{1+s}.
\end {equation}
%Moreover since $\lambda _nW_n\ll 1/n^k$ and for $s\geq -1+\delta $
%\[ W_n^s\leq \prod _{p|n}\frac {p^{(t-1)k}}{p^{(t-1)k\delta }}
%\leq n^k\prod _{p|n}\frac {1}{p^kp^{(t-1)k\delta }}\]
%we have
%\begin {equation}\label {bded2}
%\sum _{n=1}^\infty \left |\lambda _nW_n^{1+s}\right |\leq \sum _{n=1}^\infty \left (
%\prod _{p|n}\frac {1}{p^kp^{(t-1)k\delta }}\right )
%=\prod _p\left (1+\frac {1}{p^k}\sum _{t\geq 1}\frac {1}{p^{(t-1)k\Delta }}\right ).
%\end {equation}
Just as \eqref {now} is true so is
\begin {equation}\label {fodandu}
\sum _{n=1}^\infty \left |f_s^*(n)\right |=\prod _p\left (1-\frac {2}{p^{k}}\sum _{t\geq 1}
\left |\left (\frac {-(q,p^k)^s}{p^{k(1+s)}}\right )^{t-1}\right |\right ).
\end {equation}
For $-1+\delta \leq \mathfrak R\mathfrak e(s)\leq 0$ the $t$ sum here is from \eqref {hwndio} 
\[ \ll \sum _{t\geq 1}\left (\frac {1}{p^{k\delta }}\right )^{t-1}=\frac {1}{1-p^{k\delta }}\ll 1\]
so the Euler product in \eqref {fodandu} is uniformly bounded in this range and therefore
\begin {equation}\label {bded2}
\sum _{n=1}^\infty \left |f^*(n)\right |\ll 1,\hspace {5mm}\text { for }-1+\delta \leq \mathfrak R\mathfrak e(s)\leq 0.
\end {equation}
%From \eqref {hwndio2}
We have for $\mathfrak R\mathfrak e(s)\geq -1+\delta $
%(B) Suppose $\mathfrak R\mathfrak e(s)\geq -1+1/2k$.  For each prime $p$ write $A=1-2/p^k$.  As in the proof of Lemma \ref {ds}
%\[ \left |Ap^{k(1+s)}\right |\geq 1.01\]
%and therefore
%\[ \frac {1}{A+1/p^{k(1+s)}}=\frac {1}{A(1+1/Ap^{k(1+s)})}=\frac {1}{A}\sum _{l\geq 0}\left (\frac {-1}{Ap^{k(1+s)}}\right )^l
%=\sum _{l\geq 0}\frac {C_p(l)}{p^{k(1+s)l}}\]
\begin {equation}\label {11}
\frac {1}{1+1/p^{k(1+s)}}=\sum _{l\geq 0}\left (\frac {-1}{p^{k(1+s)}}\right )^l
=\sum _{l\geq 0}\frac {C_p(l)}{p^{k(1+s)l}}
\end {equation}
for some $C_p(l)$ with
\begin {equation}\label {22}
\sum _{l\geq 0}\left |\frac {C_p(l)}{p^{k(1+s)l}}\right |\ll \sum _{l\geq 0}\left (\frac {1}{p^{k\delta }}\right )^l\ll 1.
\end {equation}
as well as
\begin {equation}\label {33}
1+\frac {(q,p^k)^s}{p^{k(1+s)}}=\sum _{l'\geq 0}\frac {C_p'(l')(q,p^k)^{sl'}}{p^{k(1+s)l'}}
\end {equation}
for some $C_p'(l')$ with
\begin {equation}\label {44}
\sum _{l'\geq 0}\left |\frac {C_p'(l')(q,p^k)^{sl'}}{p^{k(1+s)l'}}\right |\leq 1+1.
\end {equation}
%\begin {equation}\label {33}
%\sum _{l'\geq 0}\left |\frac {C_p'(l')(q,p^k)^{sl'}}{p^{k(1+s)l'}}\right |\ll \left \{ \begin {array}{ll}
%1+p^{-k(1+s)}&\text { if }\mathfrak R\mathfrak e(s)\leq 0
%\\ 1+p^{-k}&\text { if }\mathfrak R\mathfrak e(s)\geq 0\end {array}\right .\hspace {3mm}\ll 1.
%\end {equation}
from \eqref {hwndio}.  From \eqref {11}, \eqref {22}, \eqref {33} and \eqref {44} there 
are for each prime $p$ and $l,l'\in \mathbb N$ some $C_p(l),C_p'(l')$ for which 
\[ \frac {1+(q,p^k)^s/p^{k(1+s)}}{1+1/p^{k(1+s)}}=
\sum _{l,l'\geq 0}\frac {C_p(l)C_p'(l')(q,p^k)^{sl'}}{p^{k(1+s)(l+l')}}\]
and
\[ \sum _{l,l'\geq 0}\left |\frac {C_p(l)C_p'(l')(q,p^k)^{sl'}}{p^{k(1+s)(l+l')}}\right |\ll 1.\]
Consequently
\begin {eqnarray*}
\prod _{p|q}\frac {1+(q,p^k)^s/p^{k(1+s)}}{1+1/p^{k(1+s)}}&=&\sum _{l_1,...,l_\omega \geq 0\atop {l_1',...,l_\omega '\geq 0}}
\frac {C_{p_1}(l_1)C_{p_1}'(l_1')\cdot \cdot \cdot C_{p_\omega }(l_\omega )C_{p_\omega }'(l_\omega ')
(q,p_1^k)^{sl_1'}\cdot \cdot \cdot (q,p_\omega ^k)^{sl_\omega '}}
{p_1^{k(1+s)(l_1+l_1')}\cdot \cdot \cdot p_\omega ^{k(1+s)(l_\omega +l_\omega ')}}
\end {eqnarray*}
and, for some $A>0$,
\[ \sum _{l_1,...,l_\omega \geq 0\atop {l_1',...,l_\omega '\geq 0}}\left |
\frac {C_{p_1}(l_1)C_{p_1}'(l_1')\cdot \cdot \cdot C_{p_\omega }(l_\omega )C_{p_\omega }'(l_\omega ')(q,p_1^k)^{sl_1'}\cdot \cdot \cdot (q,p_\omega ^k)^{sl_\omega '}}
{p_1^{k(1+s)(l_1+l_1')}\cdot \cdot \cdot p_\omega ^{k(1+s)(l_\omega +l_\omega ')}}\right |\leq A^\omega \ll \log (q+1)\]
for $\mathfrak R\mathfrak e(s)\geq -1+\delta $.  If we now define 
\[ C^*_{\mathbf l,\mathbf l'}=\prod _{i=1}^\omega C_{p_i}(l_i)C_{p_i}'(l_i'),\hspace {10mm}
W_{\mathbf l,\mathbf l'}=\left (\prod _{i=1}^\omega p_i^{l_i+l_i'}\right )^k,\hspace {10mm}
C_{\mathbf l,\mathbf l'}=\frac {C^*_{\mathbf l,\mathbf l'}}{W_{\mathbf l,\mathbf l'}},\]
\[ D_{\mathbf l'}=\prod _{i=1}^\omega (q,p_i^k)^{l_i'},\hspace {5mm}\text {and}\hspace {5mm}
Z_{\mathbf l,\mathbf l'}=\frac {D_{\mathbf l'}}{W_{\mathbf l,\mathbf l'}}\]
then 
\[ \prod _{p|q}\frac {1+(q,p^k)^s/p^{k(1+s)}}{1+1/p^{k(1+s)}}=\sum _{l_1,...,l_\omega \geq 0\atop {l_1',...,l_\omega '\geq 0}} 
C_{\mathbf l,\mathbf l'}Z_{\mathbf l,\mathbf l'}^s\]
with 
\begin {eqnarray}\label {bded}\sum _{l_1,...,l_\omega \geq 0\atop {l_1',...,l_\omega '\geq 0}}
|C_{\mathbf l,\mathbf l'}Z_{\mathbf l,\mathbf l'}^s|\ll \log (q+1)\hspace {10mm}
\end {eqnarray}
for $\mathfrak R\mathfrak e(s)\geq -1+\delta $.  The first claim now follows from \eqref {bisged} and 
the boundedness claim from \eqref {bded2} and \eqref {bded}.
\end {proof}
\newtheorem {identities}[ds]{Lemma}
\begin {identities}\label {identities}
(A)  Define
\[ \alpha =\prod _{p|q}\frac {1+(q,p^k)/p^{2k}-2/p^k}{1+1/p^{2k}-2/p^k}
\prod _p\left (1-\frac {2}{p^k}+\frac {1}{p^{2k}}\right ),\]
\[ \beta =\prod _p\left (1-\frac {1}{p^k}\right )\]
and
\[ \gamma =\prod _{p|q}\frac {1+(q,p^k)^{-1+1/k}/p-2/p^k}{1+1/p-2/p^k}
\prod _p\frac {1-2/(p^k+p^{k-1})}{1-p^{1-1/k}},\]
and let $\mathcal F^*(s)$ be as given at the start of this section.  Then
\[ \frac {\zeta (2k)\mathcal F^*(1)}{\zeta (4k)}=\alpha ,\hspace {5mm}\frac {\zeta (k)\mathcal F^*(0)}{\zeta (2k)}=\beta 
\hspace {5mm}\text {and}\hspace {5mm}\zeta (-1+1/k)\mathcal F^*(-1+1/k)=\gamma .\]
(B)  Define $\eta (q,a)$ as in \eqref {ramonramon}.  For any $q,n\in \mathbb N$
\[ \eta (q,n)=\eta \left (q,(q,n)\right )\ll q^{\epsilon -1}\]
and
\[ \sum _{a=1}^q\eta (q,a)^2=\frac {\alpha }{q}.\]
\end {identities}
\begin {proof}
(A)  So long as there are no problems with zeros of denominators we have
\begin {eqnarray*}
\mathcal F^*(s)&=&\prod _{p|q}\frac {1+(q,p^k)^s/p^{k(1+s)}-2/p^k}{1+1/p^{k(1+s)}}
\prod _{p\not |\hspace {1mm}q}\left (1-\frac {2}{p^k(1+1/p^{k(1+s)})}\right )
\\ &=&\prod _{p|q}\frac {1+(q,p^k)^s/p^{k(1+s)}-2/p^k}{1+1/p^{k(1+s)}-2/p^k}
\prod _p\left (1-\frac {2}{p^k(1+1/p^{k(1+s)})}\right ).
\end {eqnarray*}
For $s=1,0,-1+1/k$ there are clearly no problems and therefore from the Euler product expressions for the Riemann zeta function
\begin {eqnarray*}
\frac {\zeta (2k)\mathcal F^*(1)}{\zeta (4k)}&=&\prod _{p|q}\frac {1+(q,p^k)/p^{2k}-2/p^k}{1+1/p^{2k}-2/p^k}
\prod _p\left (\frac {1-p^{-4k}}{1-p^{-2k}}\right )
\left (1-\frac {2}{p^k(1+1/p^{2k})}\right )
\\ &=&\prod _{p|q}\frac {1+(q,p^k)/p^{2k}-2/p^k}{1+1/p^{2k}-2/p^k}
\prod _p\left (1+\frac {1}{p^{2k}}-\frac {2}{p^k}\right ),
\end {eqnarray*}
\begin {eqnarray*}
\frac {\zeta (k)\mathcal F^*(0)}{\zeta (2k)}&=&\prod _{p|q}\frac {1+1/p^{k}-2/p^k}{1+1/p^{k}-2/p^k}
\prod _p\left (\frac {1-p^{-2k}}{1-p^{-k}}\right )\left (1-\frac {2}{p^k(1+1/p^k)}\right )
\\ &=&\prod _p\left (1+\frac {1}{p^{k}}-\frac {2}{p^k}\right )
\end {eqnarray*}
and
\begin {eqnarray*}
&&\zeta (-1+1/k)\mathcal F^*(-1+1/k)
\\ &&\hspace {10mm}=\hspace {4mm}\prod _{p|q}\frac {1+(q,p^k)^{-1+1/k}/p-2/p^k}{1+1/p-2/p^k}
\prod _p\left (1-p^{1-1/k}\right )^{-1}\left (1-\frac {2}{p^k(1+1/p)}\right ).
\end {eqnarray*}
(B)  From \eqref {ramonramon}
\begin {eqnarray*}
\eta (q,a)&=&\sum _{D|q\atop {D\text { is ($k+1$)-free}}}\sum _{d=1\atop {(q,d^k)|a\atop {(q,d^k)=D}}}^\infty \frac {\mu (d)}{[q,d^k]}.
\end {eqnarray*}
Writing $l_0$ for the squarefree part of $D$ the $d$ sum must be
\begin {eqnarray*}
\frac {D}{q}\sum _{d=1\atop {(q,d^k)|a\atop {(q,d^k)=D}}}^\infty \frac {\mu (d)}{d^k}=
\frac {D}{q}\sum _{d=1\atop {(q,(dl_0)^k)|a\atop {(q,(dl_0)^k)=D}}}^\infty \frac {\mu (dl_0)}{(dl_0)^k}\ll \frac {D}{ql_0^k}
\end {eqnarray*}
so that 
\[ \eta (q,a)\ll \sum _{D|q\atop {D\text { is ($k+1$)-free}}}\frac {D}{ql_0^k}\ll q^{\epsilon -1}\]
which is the second claim and the first is trivial.  We have
\begin {eqnarray}\label {nionyn}
\sum _{a=1}^q\eta (q,a)^2&=&\sum _{d,d'=1}^\infty \frac {\mu (d)\mu (d')}{[q,d^k][q,{d'}^k]}
\sum _{a=1\atop {(q,d^k),(q,{d'}^k)|a}}^q1\notag 
\\ &=&q\sum _{d,d'=1}^\infty \frac {\mu (d)\mu (d')}{[q,d^k][q,{d'}^k][(q,d^k),(q,{d'}^k)]}\notag 
\\ &=&\frac {1}{q}\sum _{d,d'=1}^\infty \frac {\mu (d)\mu (d')(q,d^k,{d'}^k)}{d^k{d'}^k}\notag 
\\ &=&\frac {1}{q}\sum _{N=1}^\infty \frac {1}{N^k}\sum _{dd'=N}\mu (d)\mu (d')(q,d^k,{d'}^k)\notag 
\\ &=:&\frac {1}{q}\sum _{N=1}^\infty \frac {b_q(N)}{N^k}.
\end {eqnarray}
Clearly $b_q(N)$ is multiplicative and simple calculations show
\begin {eqnarray*}
b_q(p)=-2, 
\end {eqnarray*}
\begin {eqnarray*}
b_q(p^2)=(q,p^k) 
\end {eqnarray*}
and $b_q(p^t)=0$ for $t\geq 3$.  Consequently
\begin {eqnarray*}
\sum _{N=1}^\infty \frac {b_q(N)}{N^k}&=&\prod _p\left (1-\frac {2}{p^k}+\frac {(q,p^k)}{p^{2k}}\right )
\\ &=&\prod _{p|q}\frac {1-2/p^k+(q,p^k)/p^{2k}}{1-2/p^k+1/p^{2k}}\prod _p\left (1-\frac {2}{p^k}+\frac {1}{p^{2k}}\right )
\\ &=&\alpha 
\end {eqnarray*}
which with \eqref {nionyn} is the third claim.
\end {proof}
\newtheorem {coron}[ds]{Lemma}
\begin {coron}\label {coron}
Let $c>1$, let
\[ \mathcal A(s)=\sum _{n=1}^\infty \frac {a_n}{n^s}\]
be absolutely convergent for $\mathfrak R\mathfrak e(s)>c$, and let 
\[ A(Q)=\max _{Q/2\leq n\leq 3Q/2}|a_n|.\]
Then for $T>1$ and non-integer $Q>0$
\[ \sum _{n\leq Q}a_n\Big (Q-n\Big )=\frac {1}{2\pi i}\int _{c\pm iT}\frac {\mathcal A(s)Q^{s+1}ds}{s(s+1)}
+\mathcal O\left (\frac {QA(Q)2^c}{T}
\left (1+\frac {Q\log Q}{T}\right )+\left (1+\frac {Q^{c+1+\epsilon }}{T^2}\right )\sum _{n=1}^\infty \frac {|a_n|}{n^c}\right ).\]
In particular if $c-1\gg 1/\log Q$ then
\[ \sum _{n\leq Q}\Big (Q-n\Big )=\frac {1}{2\pi i}\int _{c\pm iT}\frac {\zeta (s)Q^{s+1}ds}{s(s+1)}
+\mathcal O\left (Q^\epsilon \left (1+\frac {Q^2}{T^2}\right )\right ).\]
\end {coron}
\begin {proof}
Take $X>0$ and define 
\[ \delta (X)=\left \{ \begin {array}{ll}0&\text { if }0<X<1\\ X-1&\text { if }X>1\end {array}\right .\]
and 
\[ I_X(T)=\frac {1}{2\pi i}\int _{c\pm iT}\frac {X^{s+1}ds}{s(s+1)}.\]
We first prove
\begin {eqnarray}\label {inte}
\left |I_X(T)-\delta (X)\right |\ll \frac {X^{c+1}}{T}\min \left \{ 1,\frac {1}{T|\log X|}\right \} .
\end {eqnarray}
Suppose first $0<X<1$ so that for $\sigma >0$ we have $X^{s+1}\ll 1$.  Then for $R>c$ 
\begin {eqnarray*}
2\pi iI_x(T)&=&-\left (\int _{c+iT}^{R+iT}+\int _{R+iT}^{R-iT}+\int _{R-iT}^{c-iT}\right )\frac {X^{s+1}ds}{s(s+1)}
\\ &\ll &\frac {1}{T^2}\int _{c}^{R}X^{\sigma +1}d\sigma +\frac {1}{R^2}\int _{\pm T}dt
\\ &\ll &\frac {X^{c+1}}{T^2|\log X|}
\end {eqnarray*}
with $R\rightarrow \infty $.  Suppose now that $X>1$ so that for $\sigma \leq -1$ we have $X^{s+1}\ll 1$.  Then for $R<-1$ 
\begin {eqnarray*}
2\pi iI_x(T)&=&Res _{s=0}\left (\frac {X^{s+1}}{s(s+1)}\right )+Res _{s=-1}\left (\frac {X^{s+1}}{s(s+1)}\right )
-\left (\int _{c+iT}^{R+iT}+\int _{R+iT}^{R-iT}+\int _{R-iT}^{c-iT}\right )\frac {X^{s+1}ds}{s(s+1)} 
%\\ &=&X-1+\mathcal O\left (\frac {1}{T^2}\int _c^\infty X^{\sigma +1}d\sigma +\frac {T}{|R|^2}\right )
%\\ &=&X-1+\mathcal O\left (\frac {X^{c+1}}{T^2\log X}\right ).
\end {eqnarray*}
and bounding the integrals as above shows 
\[ I_X(T)-(X-1)\ll \frac {X^{c+1}}{T^2|\log X|}\]
so that we can conclude that the second bound in \eqref {inte} is clear; now for the first 
bound.  If $0<X<1$ and if $\mathcal C$ is the arc of the circle 
going from $c+iT$ to $c-iT$ counterclockwise (so a circle of 
radius $\sqrt {T^2+c^2}>T$, and so that $X^{s+1}\ll X^{c+1}$ on $\mathcal C$) then
\begin {eqnarray*}
2\pi iI_X(T)&=&-\int _{\mathcal C}\frac {X^{s+1}ds}{s(s+1)} 
\\ &\ll &X^{c+1}\int _{\mathcal C}\frac {ds}{|s|\cdot |s+1|}\ll \frac {X^{c+1}}{T}.
\end {eqnarray*}
If $X>1$ the remaining part of the circle should be taken as the contour 
so that $X^{s+1}\ll X^{c+1}$ holds on the contour, and this gives a similar result.  We conclude that 
the first bound in \eqref {inte} also holds and so the proof of \eqref {inte} is complete.  Therefore by absolute convergence
\begin {eqnarray}\label {QQ}
\int _{c\pm iT}\frac {\mathcal A(s)Q^{s+1}ds}{s(s+1)}&=&
\sum _{n=1}^\infty a_nn\int _{c\pm iT}\frac {1}{s(s+1)}\left (\frac {Q}{n}\right )^{s+1}ds\notag 
\\ &=&\sum _{n=1}^\infty a_nn\delta (Q/n)+\mathcal O\left (\frac {Q^{c+1}}{T}
\sum _{n=1}^\infty \frac {|a_n|}{n^c}\min \left \{ 1,\frac {1}{T|\log (Q/n)|}\right \} \right ).\notag 
\end {eqnarray}
In general for $Z>-1$ 
\[ |\log (1+Z)|\geq \frac {|Z|}{1+Z}\]
(take logarithms of a well-known inequality to deduce $X\geq \log (1+X)$ for $X>-1$ and put in $X=-Z/(Z+1)$, which for $-1<Z\leq 0$ is 
positive and for $Z\geq 0$ satisfies $|X|\leq 1$) so that, since for $Q/2\leq n\leq 3Q/2$ we have $(n-Q)/Q>-1$,
\[ |\log (Q/n)|=\left |\log \left (1+\frac {n-Q}{Q}\right )\right |\geq \frac {|n-Q|}{n}\geq \Big \lfloor |n-Q|\Big \rfloor \Big /n.\]
Therefore
\begin {eqnarray*}
\sum _{Q/2\leq n\leq 3Q/2}\frac {|a_n|}{n^c}\min \left \{ 1,\frac {1}{T|\log (Q/n)|}\right \} &\leq &A(Q)\left (
\left (\frac {Q}{2}\right )^{-c}+\frac {2}{T}\left (\frac {Q}{2}\right )^{1-c}
\sum _{h\leq Q/2+1}\frac {1}{h}\right )
\\ &\ll &Q^{-c}A(Q)2^c\left (1+\frac {Q^{1+\epsilon }}{T}\right )\hspace {15mm}
\end {eqnarray*}
(assuming that $Q\geq 1/2$, as we can since the integral then goes into the last error 
term) and if $n$ is not in this range then $|\log (Q/n)|\gg 1$ so we deduce
\begin {eqnarray*}
\frac {Q^{c+1}}{T}\sum _{n=1}^\infty \frac {|a_n|}{n^c}\min \left \{ 1,\frac {1}{T|\log (Q/n)|}\right \} &\ll &\frac {QA(Q)2^c}{T}
\left (1+\frac {Q\log Q}{T}\right )+\frac {Q^{c+1}}{T^2}\sum _{n=1}^\infty \frac {|a_n|}{n^c}.
\end {eqnarray*}
Therefore the error term in \eqref {QQ} is of the right order of magnitude and of course the main term is
\[ \sum _{n\leq Q}a_n\Big (Q-n\Big )\]
and the main claim is proven.  For the ``in particular claim" the main claim implies an error term
\[ Q^\epsilon \left (1+\frac {Q}{T}+\left (1+\frac {Q^{c+1}}{T^2}\right )\zeta (c)\right );\]
now use $\zeta (c)\ll 1/(c-1)\ll \log Q$ and $Q^c\ll Q$.
\end {proof}
\newtheorem {ramon}[ds]{Lemma}
\begin {ramon}\label {ramon}
Take $Q>0$, $L\geq 2$ and $\Delta \in [1/2k,1/k)$.  Let 
\[ R_1=-1+\Delta \hspace {5mm}\text {and}\hspace {5mm}
R_2=\Delta k.\]
Then 
%\[ \int _{1}^{Q^3}\frac {\zeta (R_1\pm it)\zeta (R_2\pm it)Q^{\pm it}dt}{t^2}\ll \log (2+Q).\]
\[ \int _{1}^{L}\frac {\zeta (R_1+it)\zeta (R_2+it)Q^{it}dt}{t^2}\ll L^{1/4-1/2k}\log L.\]
\end {ramon}
\begin {proof}
Take $s=\sigma +it\in \mathbb C$ with $t\geq 1$ and take two parameters $N,M\gg 1$ with $NM=t/2\pi $.  Let
\[ \chi (s)=\frac {2^{s-1}\pi ^s\sec (s\pi /2)}{\Gamma (s)}.\]
By formula (4.12.3) of \cite {titchmarsh} (the definition of $\chi (s)$ comes just before) we have for $-1\leq \sigma \leq 1$
\begin {eqnarray}\label {stirling}
\chi (s)&=&\left (\frac {t}{2\pi }\right )^{1/2-\sigma -it}e^{i(t+\pi /4)}\left (1+\mathcal O\left (\frac {1}{t}\right )\right )\notag 
\\ &=&\left (\frac {t}{2\pi }\right )^{1/2-\sigma -it}e^{i(t+\pi /4)}+\mathcal O\left (\frac {1}{t^{1/2+\sigma }}%+\frac {1}{t}
\right )
\end {eqnarray}
so that 
\begin {eqnarray*}
\chi (R_2+it)\sum _{n\leq M}\frac {1}{n^{1-R_2-it}}
&=&\left (\frac {t}{2\pi }\right )^{1/2-R_2-it}e^{i(t+\pi /4)}\sum _{n\leq M}\frac {1}{n^{1-R_2-it}}
+\mathcal O\left (\frac {M^{R_2}}{t^{1/2+R_2}}\right )
\end {eqnarray*}
so by the approximate functional equation (formula (4.12.4) of \cite {titchmarsh}) 
\begin {eqnarray}\label {hanar}
\zeta (R_2+it)&=&\sum _{n\leq N}\frac {1}{n^{R_2+it}}+\chi (R_2+it)\sum _{n\leq M}\frac {1}{n^{1-R_2-it}}\notag 
\\ &&\hspace {10mm}+\hspace {4mm}\mathcal O\left (N^{-R_2}+t^{1/2-R_2}M^{R_2-1}\right )\notag 
\\ &=&\sum _{n\leq N}\frac {1}{n^{R_2+it}}+\left (\frac {t}{2\pi }\right )^{1/2-R_2-it}e^{i(t+\pi /4)}
\sum _{n\leq M}\frac {1}{n^{1-R_2-it}}\notag 
\\ &&\hspace {10mm}+\hspace {4mm}\mathcal O\left (\left (\frac {M}{t}\right )^{R_2}
\left (1+\frac {t^{1/2}}{M}\right )\right ).\hspace {15mm}
\end {eqnarray}
From the functional equation (this just preceeds 
formula (4.12.1) of \cite {titchmarsh})  and from \eqref {stirling} we have 
\begin {eqnarray*}
\zeta (R_1+it)&=&\left (\left (\frac {t}{2\pi }\right )^{1/2-R_1-it}e^{i(t+\pi /4)}
+\mathcal O\left (\frac {1}{t^{1/2+R_1}}\right )\right )\zeta (1-R_1-it)
\\ &=&\left (\frac {t}{2\pi }\right )^{1/2-R_1-it}e^{i(t+\pi /4)}\zeta (1-R_1-it)+\mathcal O\left (\frac {1}{t^{1/2+R_1}}\right )
\end {eqnarray*}
so that with \eqref {hanar} we get
\begin {eqnarray*}
&&\zeta (R_1+it)\zeta (R_2+it)
\\ &&\hspace {10mm}=\hspace {4mm}\left (\frac {t}{2\pi }\right )^{1/2-R_1-it}e^{i(t+\pi /4)}\zeta (1-R_1-it)
\sum _{n\leq N}\frac {1}{n^{R_2+it}}
\\ &&\hspace {14mm}+\hspace {4mm}\left (\frac {t}{2\pi }\right )^{1-R_1-R_2-2it}e^{2i(t+\pi /4)}\zeta (1-R_1-it)
\sum _{n\leq M}\frac {1}{n^{1-R_2-it}}
\\ &&\hspace {18mm}+\hspace {4mm}\mathcal O\left (t^{1/2-R_1}|\zeta (1-R_1-it)|\left (\frac {M}{t}\right )^{R_2}
\left (1+\frac {t^{1/2}}{M}\right )\right .
\\ &&\left .\hspace {22mm}+\hspace {4mm}
\frac {(t/M)^{1-R_2}+t^{1/2-R_2}M^{R_2}}{t^{1/2+R_1}}+\frac {1}{t^{1/2+R_1}}\left (\frac {M}{t}\right )^{R_2}
\left (1+\frac {t^{1/2}}{M}\right )\right )
\\ &&\hspace {10mm}=:\hspace {4mm}M_1(t)+M_2(t)+\mathcal O\left (t^{1/2-R_1}\left (\frac {M}{t}\right )^{R_2}
\left (1+\frac {t^{1/2}}{M}\right )\right ).
\end {eqnarray*}
Write $N=t^{1/A}$ and $M=t^{1/B}$ so the above reads
\begin {eqnarray}\label {rhyfadd}
\zeta (R_1+it)\zeta (R_2+it)=
M_1(t)+M_2(t)+\mathcal O\left (t^{1/2-R_1+R_2/B-R_2}\left (1+t^{1/2-1/B}\right )\right ).\hspace {15mm}
\end {eqnarray}
For some constant $C$
\begin {eqnarray*}
%M_1(t)x^{\pm it}&=&Ct^{1/2-R_1}\sum _{n\leq N}\sum _{m=1}^\infty \frac {e^{it(-\log t+1-\log n+\log m\pm \log x)}}
%{n^{R_2}m^{1-R_1}}
%\\ &=:&Ct^{1/2-R_1}\sum _{2\pi n^2\leq t}\sum _{m=1}^\infty \frac {e\left (f_{mx^{\pm 1}/n}(t)\right )}{n^{R_2}m^{1-R_1}}
M_1(t)Q^{it}&=&Ct^{1/2-R_1}\sum _{n\leq N}\sum _{m=1}^\infty \frac {e^{it(-\log t+1-\log n+\log m+\log Q)}}
{n^{R_2}m^{1-R_1}}
\\ &=&Ct^{1/2-R_1}\sum _{n^A\leq t\atop {2\pi nM\leq t}}\sum _{m=1}^\infty \frac {e\left (f_{mQ/n}(t)\right )}{n^{R_2}m^{1-R_1}}
\end {eqnarray*}
where 
\[ f_{X}(t)=\frac {t(-\log t+1+\log X)}{2\pi }\]
and the two summation conditions on $n$ are equivalent.  So for any $T\geq 1$
\begin {eqnarray}\label {glaw}
\int _T^{2T}\frac {M_1(t)Q^{it}dt}{t^2}&=&C\sum _{m=1}^\infty \frac {1}{m^{1-R_1}}\sum _{n^A\leq T}
\frac {1}{n^{R_2}}\int _{\max (2\pi nM,T)}^{2T}\frac {e\left (f_{mQ/n}(t)\right )dt}{t^{3/2+R_1}}.\hspace {15mm}
\end {eqnarray}
We now bound this oscillatory integral.  We have
\begin {equation}\label {differu}
2\pi f'_X(t)=-\log t+\log X.
\end {equation}
Suppose first that $T$ is large and $0<X\ll 1$.  For $\max (2\pi nM,T)<t<2T$ we have from \eqref {differu}
\[ f_X'(t)\gg 1\]
and
\[ t^{3/2+R_1}\gg T^{3/2+R_1}\]
so from Lemma 4.3 of \cite {titchmarsh}
\begin {eqnarray}\label {bach}
\int _{\max (2\pi nM,T)}^{2T}\frac {e\left (f_{X}(t)\right )dt}{t^{3/2+R_1}}\ll \frac {1}{T^{3/2+R_1}},
\hspace {10mm}\text {if }0<X\ll 1.
\end {eqnarray}
Suppose now $X$ is large.  Since from \eqref {differu} %for $t\geq 1$ 
\begin {eqnarray*}
f'_{X}(t)&\gg &|\log (t/X)|
\\ &=&\left |\log \left (1+(t-X)/X\right )\right |
\\ &\gg &\left \{ \begin {array}{ll}|t-X|/X&\text { if }t\in (X/2,3X/2)
\\ 1&\text { if not}\end {array}\right .
\\ &\gg &\left \{ \begin {array}{ll}1/\sqrt X&\text { if }t\in (X/2,X-\sqrt X)\cup (X+\sqrt X,3X/2)
\\ 1&\text { if }t\not \in (X/2,3X/2)\end {array}\right .
%\\ &\gg &\left \{ \begin {array}{ll}1/\sqrt X&\text { if }t\in (1,X-\sqrt X)\cup (X+\sqrt X,\infty )
%\\ 1&\text { if }t\not \in (X/2,3X/2)\end {array}\right .
\end {eqnarray*}
and since for $t>T$
\begin {equation}\label {erdinger}
t^{3/2+R_1}\gg T^{3/2+R_1}
\end {equation}
we have from Lemma 4.3 of \cite {titchmarsh}
\begin {eqnarray*}
&&\int _{\max (2\pi nM,T)}^{2T}\frac {e\left (f_{X}(t)\right )dt}{t^{3/2+R_1}}
\\ &&\hspace {10mm}=\hspace {4mm}\int _{\max (2\pi nM,T)\atop {t\not \in (X-\sqrt X,X+\sqrt X)}}^{2T}
+\int _{\max (2\pi nM,T)\atop {t\in (X-\sqrt X,X+\sqrt X)}}^{2T} 
\\ &&\hspace {10mm}\ll \hspace {4mm}\left \{ \begin {array}{ll}
\sqrt X/T^{3/2+R_1}&\text { if }(T,2T)\cap (X/2,3X/2)\not =\emptyset  
\\ 1/T^{3/2+R_1}&\text { if }(T,2T)\subseteq (1,\infty )\char92 (X/2,3X/2)
%\\ \sqrt X/T^{3/2+R_1}&\text { if }(T,2T)\subseteq (X-\sqrt X,X+\sqrt X).
\end {array}\right .
\\ &&\hspace {10mm}\ll \hspace {4mm}\frac {1}{T^{1+R_1}},
\end {eqnarray*}
where we have used a trivial bound for the second integral.  Therefore from \eqref {bach} 
\[ \int _{\max (2\pi nM,T)}^{2T}\frac {e\left (f_{X}(t)\right )dt}{t^{3/2+R_1}}\ll \frac {1}{T^{1+R_1}}\]
holds in fact for all $X>0$, so we deduce from \eqref {glaw}
\begin {eqnarray}\label {coffi12}
\int _T^{2T}\frac {M_1(t)Q^{it}dt}{t^2}&\ll &\frac {1}{T^{1+R_1}}
\sum _{m=1}^\infty \frac {1}{m^{1-R_1}}\sum _{n^A\leq T}\frac {1}{n^{R_2}}\notag 
\\ &\ll &\frac {1}{T^{1+R_1}}\left (T^{1/A}\right )^{1-R_2}\notag 
\\ &\ll &T^{1/A-1-R_1-R_2/A}.
\end {eqnarray}
Similarly we have
\[ \int _T^{2T}\frac {M_2(t)Q^{it}dt}{t^2}=C\sum _{m=1}^\infty \frac {1}{m^{1-R_1}}\sum _{n^B\leq T}
\frac {1}{n^{1-R_2}}\int _{\max (nN,T)}^{2T}\frac {e\left (f_{mQ/n}(t)\right )dt}{t^{1+R_1+R_2}}\]
where the oscillatroy integral is
\[ \ll \frac {1}{T^{1/2+R_1+R_2}}\]
so that 
\begin {eqnarray}\label {coffi13}
\int _T^{2T}\frac {M_2(t)Q^{it}dt}{t^2}&\ll &\frac {1}{T^{1/2+R_1+R_2}}\sum _{m=1}^\infty \frac {1}{m^{1-R_1}}\sum _{n^B\leq T}
\frac {1}{n^{1-R_2}}\notag 
\\ &\ll &\frac {1}{T^{1/2+R_1+R_2}}\left (T^{1/B}\right )^{R_2}\notag 
\\ &\ll &T^{R_2/B-1/2-R_1-R_2}.
\end {eqnarray}
Note that
\begin {equation}\label {tatws}
-1/2-R_1-R_2/2=1/2-\Delta -\Delta k/2\leq 1/2-\Delta -1/4
\end {equation}
so taking $A=B=2$ we see from \eqref {coffi12} and \eqref {coffi13}
\[ \int _T^{2T}\frac {\left (M_1(t)+M_2(t)\right )Q^{it}dt}{t^2}\ll T^{-1/2-R_1-R_2/2}\ll T^{1/4-\Delta }.\]
We assumed that $T$ is large but the bound is trivial for $T$ not large so we conclude
\begin {eqnarray*}
\int _1^{L}\frac {\left (M_1(t)+M_2(t)\right )Q^{it}dt}{t^2}&\ll &L^{1/4-\Delta }\log L
\end {eqnarray*}
and so from \eqref {rhyfadd} and \eqref {tatws}
\begin {eqnarray*}
\int _1^{L}\frac {\zeta (R_1+it)\zeta (R_2+it)Q^{it}dt}{t^2}&\ll &L^{1/4-\Delta }\log L+
\int _1^Lt^{-3/2-R_1-R_2/2}dt
\\ &\ll &L^{1/4-\Delta }\log L.
\end {eqnarray*}
\end {proof}
\newtheorem {coro}[ds]{Lemma}
\begin {coro}\label {coro}
Let $\alpha ,\beta ,$ and $\gamma $ be as in Lemma \ref {identities} and let $\mathcal F^*(s)$ be as given at the start of this 
section.  For $X>0$ and $T,c>1$ 
\begin {eqnarray*}
&&\int _{c\pm iT}\frac {\zeta (s)\zeta \left (k(s+1)\right )\mathcal F^*(s)X^{s+1}ds}
{s(s+1)\zeta \left (2k(s+1)\right )}
\\ &&\hspace {10mm}=\hspace {4mm}\frac {\alpha X^2}{2}-\frac {\beta X}{2}
+\frac {k\gamma X^{1/k}}{(-1+1/k)\zeta (2)}
\\ &&\hspace {14mm}+\hspace {4mm}
\mathcal O\left (T^\epsilon \left (q^\epsilon T^{1/4}\left (\frac {X}{T}\right )^{1/2k}+\frac {X^{c+1}}{T^2}+1\right )\right ).
\end {eqnarray*}
\end {coro}
\begin {proof}
%An obvious bound for the integrand is $X^{c+1}/t^{3/2}$ so the integral is $\ll X^{c+1}$ and therefore we 
%can assume $T$ and $X$ are large.  
For $s\in \mathbb C$ write always $s=\sigma +it$ for $\sigma ,t\in \mathbb R$ and let
\begin {equation}\label {I}
\mathcal I(s)=\frac {\zeta (s)\zeta \left (k(s+1)\right )\mathcal F^*(s)}{\zeta \left (2k(s+1)\right )}.
\end {equation}
%For $\sigma =2$ the functions $\Delta (s),\mathcal F^*(s)$ and the zeta 
%factors are all absolutely bounded and so the same is true of $\mathcal I(s)$, so for future reference we record 
%\begin {eqnarray}\label {tailend}
%\int _{2\pm i\infty \atop {|t|\geq T}}\frac {\mathcal I(s)X^{s+1}ds}{s(s+1)}\ll X^3\int _{T}^\infty \frac {dt}{t^2}\ll 1.
%\end {eqnarray}
Let $R_1=-1+1/2k+\tau $ for some $0<\tau <1/k$.  We have already established (just 
before Lemma \ref {ds}) that $\mathcal F^*(s)\ll q^\epsilon $ for $\sigma \geq -1+\delta $, therefore
\begin {equation}\label {f}
\mathcal I(s)\ll \frac {q^\epsilon |\zeta (s)\zeta \left (k(s+1)\right )|}
{|\zeta \left (2k(s+1)\right )|},\hspace {5mm}\text { for }\sigma \geq R_1.
\end {equation}
%In particular for $\sigma \geq 0$ we have $\mathcal I(s)\ll |\zeta (s)|\ll |t|^{1/2+\epsilon }$ so 
%that, since $\Delta (1),\Delta (0),\Delta (-1+1/k)\ll 1$, the claim is trivial for $X\ll 1$ so we may assume $X$ is large.
%\\
%\\
On $\mathfrak R\mathfrak e(s)\geq -1+\delta $ we know by Lemma \ref {ds} that $\mathcal I(s)$ is holomorphic except 
for simple poles at $s=1$ and $s=-1+1/k$ so by the Residue Theorem 
\begin {eqnarray}\label {integral}
2\pi i\int _{c\pm iT}\frac {\mathcal I(s)X^{s+1}ds}{s(s+1)}&=&
\frac {X^2Res_{s=1}\mathcal I(s)}{2}+\mathcal I(0)X
+\frac {kX^{1/k}Res_{s=-1+1/k}\mathcal I(s)}{-1+1/k}\notag 
\\ &&-\hspace {4mm}2\pi i\left (\int _{c+iT}^{R_1+iT}+\int _{R_1+iT}^{R_1-iT}
+\int _{R_1-iT}^{c-iT}\right )\frac {\mathcal I(s)X^{s+1}ds}{s(s+1)}.\hspace {15mm}
\end {eqnarray}
%here obviously $\mathcal I(1)$ means the value of the 
%residue of $\mathcal I(s)$ at $s=1$ and similarly for $\mathcal I(-1+1/k)$.  
It is standard that for $t\geq 1$
\[ \zeta (s)\ll t^\epsilon \left \{ \begin {array}{ll}t^{1/2-\sigma }&\text { for }\sigma \leq 0
\\ \max \{ 1,t^{1/2-\sigma /2}\} &\text { for }\sigma \geq 0
\\ t^{1/4}&\text { for }\sigma \geq 1/2\end {array}\right .\]
and
\[ \zeta (\sigma )\ll \left \{ \begin {array}{ll}1&\text { for }\sigma \geq 2k
\\ 1/|\sigma -1|&\text { for }1\leq \sigma \leq 2;\end {array}\right .\]
we will now use these bounds freely without comment.  If $0\leq \sigma \leq 2$ and $t\geq 1$ we have 
\[ \zeta (s)\ll t^\epsilon \max \{ 1,t^{1/2-\sigma /2}\} ,\]
\[ \zeta \left (k(s+1)\right )\ll 1\]
and
\[ \frac {1}{\zeta \left (2k(s+1)\right )}\ll \zeta \left (2k(\sigma +1)\right )\ll 1,\]
so from \eqref {f}
\[ \mathcal I(s)\ll t^\epsilon \max \{ 1,t^{1/2-\sigma /2}\} \]
and therefore
\begin {equation}\label {h1}
\int _{iT}^{c+iT}\frac {\mathcal I(s)X^{s+1}ds}{s(s+1)}\ll T^\epsilon \left (\frac {X}{T^{3/2}}+\frac {X^{c+1}}{T^2}\right ).
\end {equation}
If $R_1\leq \sigma \leq 0$ then for $t\geq 1$
\[ \zeta (s)\ll t^{1/2-\sigma },\]
\[ \zeta \left (k(s+1)\right )\ll t^{1/2}\]
and
\[ \frac {1}{\zeta \left (2k(s+1)\right )}\ll \zeta \left (2k(\sigma +1)\right )\ll \frac {1}{|2k(\sigma +1)-1|}
\ll \frac {1}{\tau },\]
so from \eqref {f}
\begin {equation}\label {garlic}
\mathcal I(s)\ll \frac {t^{1-\sigma }}{\tau }
\end {equation}
and therefore
\begin {equation}\label {h2}
\int _{R_1+iT}^{iT}\frac {\mathcal I(s)X^{s+1}ds}{s(s+1)}\ll \frac {1}{\tau }\left (\frac {X^{R_1+1}}{T^{1+R_1}}+\frac {X}{T}\right )
\ll \frac {1}{\tau }\left (1+\frac {X}{T}\right ).
\end {equation}
From \eqref {h1} and \eqref {h2} we have
\begin {equation}\label {inte1}
\left (\int _{c+iT}^{R_1+iT}+\int _{R_1-iT}^{c-iT}\right )\frac {\mathcal I(s)X^{s+1}ds}{s(s+1)}
\ll \frac {1}{\tau }\left (1+\frac {X}{T}+\frac {X^{c+1}}{T^2}\right )
\ll \frac {T^\epsilon }{\tau }\left (1+\frac {X^{c+1}}{T^2}\right )
\end {equation}
a similar argument for the second integral obviously valid.  We now turn to the vertical 
contribution in \eqref {integral}.  Denote by $\omega $ the number of prime factors of $q$.  For 
given integers $n,l_1,...,l_{\omega },l_1',...,l_\omega '\geq 0$ 
write $\mathbf n=(n,l_1,...,l_{\omega },l_1',...,l_\omega ')$.  Let $W_n,Z_{\mathbf l,\mathbf l'}$ be as in 
Lemma \ref {ds2}.  Then that lemma says that for given $\mathbf n$ 
there are $a_\mathbf n=a_\mathbf n(\sigma )\in \mathbb R$ such that for $-1+\delta \leq \sigma \leq 0$
\[ \mathcal F^*(s)=\sum _{\mathbf n}a_\mathbf n\left (W_nZ_{\mathbf l,\mathbf l'}\right )^{it}\]
and
\begin {equation}\label {doch}
\sum _{\mathbf n}|a_\mathbf n|\ll 1.
\end {equation}
Therefore
\begin {eqnarray*}
\frac {\mathcal F^*(R_1+it)X^{it}}{\zeta \left (2k(R_1+it+1)\right )}&=&
\sum _{m,\mathbf n}\frac {\mu (m)a_\mathbf n}{m^{2k(R_1+1)}}\left (\frac {XW_nZ_{\mathbf l,\mathbf l'}}{m^{2k}}\right )^{it},
\end {eqnarray*}
so from \eqref {I}, Lemma \ref {ramon} and \eqref {doch}
\begin {eqnarray}\label {bread}
%&&\int _1^T\frac {\mathcal I(R_1+it)X^{it}dt}{t^2}\notag 
%\\ &&\hspace {10mm}=\hspace {4mm}\sum _{n,m\geq 1\atop {l_1,...,l_\omega \geq 0\atop {l_1',...,l_\omega '\geq 0}}}
%\frac {C_{\mathbf l,\mathbf l'}Z_{\mathbf l,\mathbf l'}(q)^{R_1}f_{R_1+1}(n)}{W_{\mathbf l,\mathbf l'}(q)^{k}m^{2k(R_1+1)}}
%\int _1^T\frac {\zeta (R_1+it)\zeta \left (k(R_1+1+it)\right )}{t^2}
%\Big (XZ_{\mathbf l,\mathbf l'}(q)\Big )^{it}dt\notag 
%\\ &&\hspace {10mm}\ll \hspace {4mm}\sum _{n,m\geq 1\atop {l_1,...,l_\omega \geq 0\atop {l_1',...,l_\omega '\geq 0}}}
%\frac {C_{\mathbf l,\mathbf l'}Z_{\mathbf l,\mathbf l'}(q)^{R_1}f_{R_1+1}(n)}{W_{\mathbf l,\mathbf l'}(q)^{k}m^{2k(R_1+1)}}
%\log \left (2+XZ_{\mathbf l,\mathbf l'}(q)\right )\notag 
%\\ &&\hspace {10mm}\ll \hspace {4mm}\log \left (2+X\right ),
\int _1^T\frac {\mathcal I(R_1+it)X^{it}dt}{t^2}&=&\sum _{m,\mathbf n}\frac {\mu (m)a_\mathbf n}{m^{2k(R_1+1)}}
\int _1^T\frac {\zeta (R_1+it)\zeta \left (k(R_1+it+1)\right )}{t^2}
\left (\frac {XW_nZ_{\mathbf l,\mathbf l'}}{m^{2k}}\right )^{it}dt\notag 
\\ &\ll &T^{1/4-1/2k}\log T\sum _{m,\mathbf n}\left |\frac {\mu (m)a_\mathbf n}{m^{2k(R_1+1)}}\right |\notag 
\\ &\ll &T^{1/4-1/2k}(\log T)\zeta (1+2k\tau )\ll \frac {T^{1/4-1/2k}\log T}{\tau }.
\end {eqnarray}
We clearly have for $\sigma \geq -1+\delta $
\[ \mathcal I(s)\ll q^\epsilon \left \{ \begin {array}{ll}1&\text { for }0\leq t\leq 1
\\ t^{7/4}&\text { for }t\geq 1\end {array}\right .\]
and for $t\geq 1$ we have
\[ \frac {1}{s(s+1)}=\frac {1}{t^2}+\mathcal O\left (\frac {1}{t^3}\right ),\]
therefore from \eqref {bread} 
\begin {eqnarray*}
\int _{R_1}^{R_1+iT}\frac {\mathcal I(s)X^{s+1}ds}{s(s+1)}&=&X^{R_1+1}\int _{1}^{T}\frac {\mathcal I(R_1+it)X^{it}ds}{t^2}\notag 
\\ &&\hspace {4mm}+\hspace {4mm}\mathcal O\left (X^{R_1+1}\int _{R_1\atop {t\geq 1}}^{R_1+i\infty }
\frac {|\mathcal I(s)|ds}{t^3}+X^{R_1+1}\int _{R_1}^{R_1+i}\frac {|\mathcal I(s)|ds}{|s(s+1)|}\right )\notag 
\\ &\ll &\frac {X^{R_1+1}T^{1/4-1/2k}\log T}{\tau }+q^\epsilon X^{R_1+1}
\\ &=&\frac {X^\tau T^{1/4+\epsilon }}{\tau }\left (\frac {X}{T}\right )^{1/2k}q^\epsilon .
\end {eqnarray*}
A similar bound obviously holding 
for $t$ negative we conclude
\begin {equation}\label {inte2}
\int _{R_1+iT}^{R_1-iT}\frac {\mathcal I(s)X^{s+1}ds}{s(s+1)}
\ll \frac {X^\tau T^{1/4+\epsilon }}{\tau }\left (\frac {X}{T}\right )^{1/2k}q^\epsilon .
\end {equation}
From Lemma \ref {identities} (A) we have
\[ Res _{s=1}\mathcal I(s)=\frac {\zeta (2k)\mathcal F^*(1)}{\zeta (4k)}=\alpha ,\]
\[ \mathcal I(0)=\frac {\zeta (0)\zeta (k)\mathcal F^*(0)}{\zeta (2k)}=-\frac {\beta }{2}\]
and
\[ Res _{s=-1+1/k}\mathcal I(s)=\frac {\zeta (-1+1/k)\mathcal F^*(-1+1/k)}{\zeta (2)}
=\frac {\gamma }{\zeta (2)}\]
so the main terms in \eqref {integral} are
\[ \frac {\alpha X^2}{2}-\frac {\beta X}{2}+\frac {k\gamma X^{1/k}}{(-1+1/k)\zeta (2)}=:M(X).\]
This with \eqref {inte1} and \eqref {inte2} means \eqref {integral} becomes
\begin {eqnarray*}
\int _{c\pm iT}\frac {\mathcal I(s)X^{s+1}ds}{s(s+1)}&=&
M(X)+\mathcal O\left (\frac {T^\epsilon }{\tau }\left (q^\epsilon X^\tau T^{1/4}\left (\frac {X}{T}\right )^{1/2k}
+1+\frac {X^{c+1}}{T^2}\right )\right )
\\ &=&M(X)+\mathcal O\left (T^\epsilon \left (q^\epsilon T^{1/4}
\left (\frac {X}{T}\right )^{1/2k}+1+\frac {X^{c+1}}{T^2}\right )\right )
\end {eqnarray*}
on taking $\tau =1/\log X$, so long as $X$ is large.  If $X$ is not large then the claim is trivial, the integrand being 
trivially $\ll t^{\epsilon -2}$ for $\sigma =c$.
\end {proof}
\newtheorem {fine}[ds]{Lemma}
\begin {fine}\label {fine}
For any $x,y>0$
\[ \sum _{[d,d']\leq y}1\leq y^{1+\epsilon },\]
\[ \sum _{[d,d']>y}\frac {1}{[d^k,{d'}^k]}\ll y^{1-k+\epsilon }\]
and, for $N\leq x$,
\[ \sum _{d,d\atop {[d,d']>y}}\sum _{n\leq x\atop {n\equiv 0(d^k)\atop {n\equiv -N({d'}^k)}}}1
\ll xy^{1-k+\epsilon }+x^{2/(k+1)+\epsilon }.\] 
\end {fine}
\begin {proof}
Since 
\[ \sum _{[d,d']=n}1\ll n^\epsilon \]
we have
\[ \sum _{[d,d']\leq y}1=\sum _{n\leq y}\sum _{[d,d']=n}1\ll y^{1+\epsilon }\]
and
\[ \sum _{[d,d']>y}\frac {1}{[d^k,{d'}^k]}=\sum _{n>y}\frac {1}{n^k}\sum _{[d,d']=n}1\ll y^{1-k+\epsilon }\]
which are the first two claims.  Let $Z$ be a parameter.  We have with a divisor estimate
\begin {eqnarray*}
\sum _{d,d'\atop {d>Z}}\sum _{n\leq x\atop {n\equiv 0(d^k)\atop {n\equiv -N({d'}^k)}}}1
\ll x^\epsilon \sum _{d^k\leq x\atop {d>Z}}\sum _{n\leq x\atop {n\equiv 0(d^k)}}1
\ll x^{1+\epsilon }\sum _{d>Z}\frac {1}{d^k}\ll x^{1+\epsilon }Z^{1-k}
\end {eqnarray*}
and similarly for the terms with $d'>Z$.  On the other hand the second claim implies
\begin {eqnarray*}
\sum _{d,d'\leq Z\atop {[d,d']>y}}\sum _{n\leq x\atop {n\equiv 0(d^k)\atop {n\equiv -N({d'}^k)}}}1&\ll &
\sum _{d,d'\leq Z\atop {[d,d']>y}}\left (\frac {x}{[d^k,{d'}^k]}+1\right )
\\ &\ll &xy^{1-k+\epsilon }+Z^2
\end {eqnarray*}
%On the other hand 
%\[ \sum _{d'>(y/d)^{1/k}}\frac {(d,d')^k}{{d'}^k}=\sum _{n|d}n^k\sum _{d'>(y/d)^{1/k}\atop {(d,d')=n}}\frac {1}{{d'}^k}
%\leq \sum _{n|d}\sum _{nd'>(y/d)^{1/k}}\frac {1}{{d'}^k}\ll y^{1/k-1}d^{k-1+\epsilon }\]
%so that
%\begin {eqnarray*}
%\sum _{d,d'\leq Z\atop {dd'>y}}\sum _{n\leq x\atop {n\equiv 0(d^k)\atop {n\equiv -ql({d'}^k)}}}1&\ll &
%\sum _{d,d'\leq Z\atop {dd'>y}}\left (\frac {x}{[d^k,{d'}^k]}+1\right )
%\\ &\ll &x\sum _{d\leq Z}\frac {1}{d^k}\sum _{d'>y/d}\frac {(d,d')^k}{{d'}^k}+Z^2
%\\ &\ll &xy^{1-k}\sum _{d\leq Z}d^{\epsilon -1}+Z^2
%\\ &\ll &xy^{1-k}Z^\epsilon +Z^2
%\end {eqnarray*}
and therefore
\begin {eqnarray*}
\sum _{[d,d']>y}\sum _{n\leq x\atop {n\equiv 0(d^k)\atop {n\equiv -N({d'}^k)}}}1&\ll &
xy^{1-k+\epsilon }+Z^2+x^{1+\epsilon }Z^{1-k}
\end {eqnarray*}
which gives the claim on choosing $Z=x^{1/(k+1)}$.
\end {proof}
\addcontentsline {toc}{subsection}{3.2 - Proof of theorem}
\section* {3.2 - Proof of theorem}
Let $1\leq q\leq x$ and define $\eta (q,a)$ and $V_x(q)$ as in \eqref {ramonramon} and \eqref {ramonramon2}.  Opening 
the square we have
\begin {eqnarray}\label {dechra}
V_x(q)&=&\sum _{a=1}^q\sum _{n,n'\leq x\atop {n,n'\in \mathcal S\atop {n\equiv n'\equiv a(q)}}}1
-2x\sum _{a=1}^q\eta (q,a)\sum _{n\leq x\atop {n\in \mathcal S\atop {n\equiv a(q)}}}1+x^2\sum _{a=1}^q\eta (q,a)^2\notag 
\\ &=&\sum _{n,n'\leq x\atop {n,n'\in \mathcal S\atop {n\equiv n'(q)}}}1
-2x\sum _{n\leq x\atop {n\in \mathcal S}}\eta (q,n)+x^2\sum _{a=1}^q\eta (q,a)^2\notag 
\\ &=:&A_x(q)-2xB_x(q)+x^2\sum _{a=1}^q\eta (q,a)^2.
\end {eqnarray}
From Lemma \ref {identities} (B) we have $\eta (q,n)=\eta \left (q,(q,n)\right )$ and $\eta (q,d)\ll 1$.  Therefore from 
Lemma 2.2 (ii) of \cite {kfree} we have for some constants $c_{dh},c_q$ and a new parameter $X\geq 1$
\begin {eqnarray}\label {B}
B_X(q)&=&\sum _{d|q}\eta (q,d)\sum _{n\leq X\atop {n\in \mathcal S\atop {(n,q)=d}}}1\notag 
\\ &=&\sum _{d|q}\eta (q,d)\sum _{h|q/d}\mu (h)\sum _{n\leq X\atop {n\in \mathcal S\atop {dh|n}}}1\notag 
\\ &=&X\sum _{d|q}\eta (q,d)\sum _{h|q/d}\mu (h)c_{dh}
+\mathcal O\left (X^{1/k+\epsilon }\sum _{d|q}|\eta (q,d)|\sum _{h|q/d}|\mu (h)|\right )\notag 
\\ &=&Xc_q+\mathcal O\left (X^{1/k+\epsilon }\right )\notag 
\\ &\sim &Xc_q
\end {eqnarray}
with $X\rightarrow \infty $.  But it is easy to establish 
\[ \sum _{n\leq X\atop {n\in \mathcal S\atop {n\equiv a(q)}}}1\sim X\eta (q,a)\]
so that evidently 
\[ B_X(q)=\sum _{a=1}^q\eta (q,a)\sum _{n\leq X\atop {n\in \mathcal S\atop {n\equiv a(q)}}}1\sim X\sum _{a=1}^q\eta (q,a)^2\]
so \eqref {B} implies 
\[ c_q=\sum _{a=1}^q\eta (q,a)^2\]
and therefore the last but one line of \eqref {B} says
\begin {equation}\label {BB}
B_x(q)=x\sum _{a=1}^q\eta (q,a)^2+\mathcal O\left (x^{1/k+\epsilon }\right ).
\end {equation}
It is well known that
\[ \sum _{n\leq x\atop {n\in \mathcal S}}1=\frac {x}{\zeta (k)}+\mathcal O\left (x^{1/k}\right )\]
therefore
\begin {eqnarray}\label {deffi}
A_x(q)&=&2\sum _{n<n'\leq x\atop {n,n'\in \mathcal S\atop {n\equiv n'(q)}}}1+\sum _{n\leq x\atop {n\in \mathcal S}}1\notag 
\\ &=&2\sum _{l\leq x/q}\sum _{n,n'\leq x\atop {n,n'\in \mathcal S\atop {n'-n=ql}}}1
+\frac {x}{\zeta (k)}+\mathcal O\left (x^{1/k}\right )\notag 
\\ &=:&2C_x(q)+\frac {x}{\zeta (k)}+\mathcal O\left (x^{1/k}\right )
\end {eqnarray}
so we deduce from \eqref {dechra} and \eqref {BB}
\begin {eqnarray}\label {va}
V_x(q)&=&2C_x(q)+\frac {x}{\zeta (k)}-x^2\sum _{a=1}^q\eta (q,a)^2+\mathcal O\left (x^{1/k+\epsilon }\right ).
\end {eqnarray}
Take a parameter $y\leq x^{1/k}$ so that $[d,d']\leq y$ is a stronger condition than $d^k,{d'}^k\leq x$.  Using
\[ \sum _{d^k|n}\mu (d)=\left \{ \begin {array}{ll}1&\text { if $n$ is $k$-free}\\ 0&\text { if not}\end {array}\right .\]
we see that
\begin {eqnarray*}
\sum _{n,n'\leq x\atop {n,n'\in \mathcal S\atop {n'-n=ql}}}1&=&\sum _{d,d'\leq x}\mu (d)\mu (d')
\sum _{n,n'\leq x\atop {n\equiv 0(d^k)\atop {n'\equiv 0(d'^k)\atop {n'-n=ql}}}}1
\\ &=&\sum _{d,d'\leq x}\mu (d)\mu (d')\sum _{n\leq x-ql\atop {n\equiv 0(d^k)\atop {n\equiv -ql(d'^k)}}}1
\\ &=&\sum _{[d,d']\leq y\atop {(d^k,d'^k)|ql}}\mu (d)\mu (d')\left (\frac {x-ql}{[d^k,d'^k]}+\mathcal O(1)\right )
+\mathcal O\left (\sum _{[d,d']>y}\sum _{n\leq x\atop {n\equiv 0(d^k)\atop {n\equiv -ql(d'^k)}}}1\right )
\\ &=&(x-ql)\sum _{d,d'=1\atop {(d^k,d'^k)|ql}}^\infty \frac {\mu (d)\mu (d')}{[d^k,d'^k]}
+\mathcal O\left (\sum _{[d,d']\leq y}1\right )
\\ &&+\hspace {4mm}\mathcal O\left ((x-ql)\sum _{[d,d']>y}\frac {1}{[d^k,d'^k]}\right )
+\mathcal O\left (\sum _{[d,d']>y}\sum _{n\leq x\atop {n\equiv 0(d^k)\atop {n\equiv -ql(d'^k)}}}1\right ).
\end {eqnarray*}
From Lemma \ref {fine} the error terms here are for $ql\leq x$
\[ y^{1+\epsilon }+xy^{1-k+\epsilon }+x^{2/(k+1)+\epsilon }\ll x^{2/(k+1)+\epsilon }\]
after setting $y=x^{1/k}$, so that 
\[ \sum _{n,n'\leq x\atop {n,n'\in \mathcal S\atop {n'-n=ql}}}1=
(x-ql)\sum _{d,d'=1\atop {(d^k,d'^k)|ql}}^\infty \frac {\mu (d)\mu (d')}{[d^k,d'^k]}
+\mathcal O\left (x^{2/(k+1)+\epsilon }\right )\]
so from \eqref {deffi}
\begin {eqnarray}\label {peter}
C_x(q)&=&\sum _{d,d'}\frac {\mu (d)\mu (d')}{[d^k,d'^k]}\sum _{l\leq x/q\atop {(d^k,d'^k)|ql}}\Big (x-ql\Big )
+\mathcal O\left (\left (x^{2/(k+1)+\epsilon }\right )\sum _{l\leq x/q}1\right )\notag 
\\ &=&\sum _{d,d'}\frac {\mu (d)\mu (d')[q,(d^k,d'^k)]}{[d^k,d'^k]}\sum _{l\leq x/[q,(d^k,d'^k)]}
\left (\frac {x}{[q,(d^k,d'^k)]}-l\right )+\mathcal O\left (\frac {x^{1+2/(k+1)+\epsilon }}{q}\right )\notag 
\\ &=:&\mathcal J(x)+\mathcal O\left (\frac {x^{1+2/(k+1)+\epsilon }}{q}\right ).
\end {eqnarray}
From now on all $\ll $ symbols will denote bounds up to $x^\epsilon $ bounds so that \eqref {va} and \eqref {peter} read
\begin {equation}\label {suncream}
V_q(x)=2\mathcal J(x)+\frac {x}{\zeta (k)}-x^2\sum _{a=1}^q\eta (q,a)^2+\mathcal O\left (\frac {x^{1+2/(k+1)}}{q}\right ).
\end {equation}
Assuming as we can that $x$ is not an integer, write $Q=x/[q,(d^k,d'^k)]$ and let $c=1+1/\log Q$.  From 
Lemma \ref {coron} the inner sum in $\mathcal J(x)$ is for any $T>1$
\[ \int _{c\pm iT}\frac {\zeta (s)}{s(s+1)}\left (\frac {x}{[q,(d^k,d'^k)]}\right )^{s+1}ds
+\mathcal O\left (1+\left (\frac {x}{T[q,(d^k,d'^k)]}\right )^{2}\right )\]
so from Lemma \ref {ds} and Lemma \ref {coro}
\begin {eqnarray*}
\mathcal J(x)&=&%\sum _{d,d'}\frac {\mu (d)\mu (d')[q,(d^k,d'^k)]}{[d^k,d'^k]}
%\int _{2\pm \infty }\frac {\zeta (s)}{s(s+1)}\left (\frac {x}{[q,(d^k,d'^k)]}\right )^{s+1}ds
%\\ &=&
\int _{c\pm iT}\frac {\zeta (s)x^{s+1}}{s(s+1)}
\left (\sum _{d,d'}\frac {\mu (d)\mu (d')}{[d^k,d'^k][q,(d^k,d'^k)]^s}\right )ds
+\mathcal O\left (\sum _{d,d'}\frac {|\mu (d)\mu (d')|}{[d^k,d'^k]}\left (1+\frac {x^2}{T^2[q,(d^k,d'^k)]}\right )\right )
\\ &=&q\int _{c\pm iT}\frac {\zeta (s)\zeta \left (k(s+1)\right )\mathcal F^*(s)}
{s(s+1)\zeta \left (2k(s+1)\right )}\left (\frac {x}{q}\right )^{s+1}ds
+\mathcal O\left (1+\frac {x^2}{qT^2}\sum _{d,d'=1}^\infty \frac {(q,d^k,{d'}^k)}{d^k{d'}^k}\right )
\\ &=&\frac {\alpha x^2}{2q}-\frac {\beta x}{2}
+\frac {k\gamma q^{1-1/k}x^{1/k}}{(-1+1/k)\zeta (2)}+
\mathcal O\left (q\left (T^{1/4}\left (\frac {x}{qT}\right )^{1/2k}+\left (\frac {x}{qT}\right )^2+1\right )\right )
\end {eqnarray*}
where $\alpha ,\beta ,\gamma $ are as in Lemma \ref {identities}, and assuming $T\leq x^2$.  Setting
\[ T=\left (\frac {x}{q}\right )^{V}\]
where
\[ V=\frac {2-1/2k}{9/4-1/2k}\]
the error term becomes 
\[ \ll q\left (\frac {x}{q}\right )^{2/(9-2/k)}\]
and so from \eqref {suncream}
\begin {eqnarray}\label {diwadd}
&&V_x(q)=\left (\frac {\alpha }{q}-\sum _{a=1}^q\eta (q,a)^2\right )x^2
+\left (\frac {1}{\zeta (k)}-\beta \right )x
+\frac {2k\gamma q^{1-1/k}x^{1/k}}{(-1+1/k)\zeta (2)}\notag 
\\ &&\hspace {20mm}+\hspace {4mm}\mathcal O\left (q\left (\frac {x}{q}\right )^{2/(9-2/k)}\right )
+\mathcal O\left (\frac {x^{1+2/(k+1)}}{q}\right ).\notag 
\end {eqnarray}
From Lemma \ref {identities} (B) the $x^2$ coefficient vanishes.   Directly from the 
definitions (Lemma \ref {identities}) we see that $\beta =\zeta (k)^{-1}$ so 
the $x$ coefficient also vanishes.  Again from the definitions the $x^{1/k}$ coefficient is
\[ \frac {2kq^{1-1/k}}{(-1+1/k)\zeta (2)}\prod _p\left (\frac {1-2/(p^k+p^{k-1})}{1-p^{1-1/k}}\right )
\prod _{p|q}\left (\frac {1+(q,p^k)^{-1+1/k}/p-2/p^k}{1+1/p-2/p^k}\right )\]
and we have our theorem.

$\hspace {1mm}$
\\
\\
\\
\\ \emph {e-mail address} - tomos.parry1729@hotmail.co.uk

\end {document}